\documentclass[12pt]{article}
\usepackage{amssymb,latexsym,amscd,amsmath,amsfonts,enumerate,supertabular}
\usepackage{graphicx}
\usepackage{tabularx}
\usepackage{caption}
\newtheorem{thm}{Theorem}[section]

 \newtheorem{prop}[thm]{Proposition}
 \newtheorem{rem}[thm]{Remark}
 \numberwithin{equation}{section}
\begin{document}
\title{{\bf A novel approach to fully third order nonlinear boundary value problems}}
 
 \author{Dang Quang A $^1$, Dang Quang Long $^2$ \\
 $^1$Center for Informatics and Computing, VAST\\ 18 Hoang Quoc Viet, Cau giay, Hanoi, Vietnam\\
E-mail: dangquanga@cic.vast.vn\\
$^2$ Institute of Information Technology, VAST\\
 E-mail: dqlong88@gmail.com} 
\date{}       
\maketitle       


\begin{abstract}
\small  
In this work we propose a novel approach to investigate boundary value problems (BVPs) for fully third order differential equations. It is based on the reduction of BVPs to operator equations for the nonlinear terms but not for the functions to be sought. By this approach we have established the existence, uniqueness, positivity and monotony of solutions and the convergence of the iterative method for approximating the solutions under some easily verified conditions in bounded domains. These conditions are much simpler and weaker than those of other authors for studying solvability of the problems before by using different methods. Many examples illustrate the obtained theoretical results.

\small  
\end{abstract}
\small 

\normalsize 

 \section {Introduction}
 In recent years the boundary value problems (BVPs) for third order nonlinear differential equations  
 have attracted attention from many researchers. A number of works are devoted to the existence, uniqueness and positivity of solutions of the problems with different boundary conditions. The methods for investigating qualitative aspects of the problems are diverse, including the method of lower and upper solutions and monotone technique \cite{Cabada94,Cabada95, Feng09,Feng-Liu,Gross,Yao-Feng}, Leray-Schauder continuation principle \cite{Hopkins}, fixed point theory on cones \cite{Li2010}, etc.. It should be emphasized that in the above works there is an essential assumption that the function $f(t,x,y,z): \; [0, 1] \times \mathbb{R}^3  \rightarrow \mathbb{R}$ satisfies a Nagumo-type condition on the last two variables \cite{Li-Li}, or linear growth in $x, y, z$ at infinity \cite{Hopkins}, or some complicated conditions including monotone increase in each of $x$ and $y$ \cite{Bai}, or one-side Lipschitz condition  in $x$ for $f=f(t,x)$ \cite{Yao-Feng} and in $x, y$ for $f=f(t,x,y)$ \cite{Feng-Liu}. Sun et al. in \cite{Sun} studied the existence of monotone positive solution of the BVP for the case  $f=f(u(t))$ under  conditions, which are difficult to be verified.\par

Motivated greatly by the above-mentioned works, in this paper we propose an efficient method to investigate the solvability and approximation of BVPs for the fully third order equation 
 \begin{equation}\label{eqfull}
 u'''(t)=f(t,u(t),u'(t),u''(t)), \quad 0<t<1
 \end{equation}
 with general boundary conditions
 \begin{equation}\label{bc1}
 \begin{split}
 B_1[u]=\alpha_1 u(0)+\beta_1 u'(0) + \gamma_1 u''(0) =0,\\
 B_2[u]=\alpha_2 u(0)+\beta_2 u'(0) + \gamma_2 u''(0) =0,\\
 B_3[u]=\alpha_1 u(1)+\beta_1 u'(1) + \gamma_1 u''(1) =0,\\
 \end{split}
 \end{equation}
 such that

 \begin{equation}\label{rank}
 Rank
\begin{pmatrix}
\alpha_1 & \beta_1 & \gamma_1 & 0 &0&0\\
\alpha_2 & \beta_2 & \gamma_2 &0 &0&0\\
0& 0 & 0 & \alpha_3& \beta_3 & \gamma_3\\
\end{pmatrix}=3.
\end{equation}
The boundary conditions \eqref{bc1} include as particular cases the boundary conditions considered in \cite{Bai,Feng09,Feng-Liu,Hopkins,Li-Li,Yao-Feng}. Notice that if there are given one boundary condition at $t=0$ and two boundary conditions at $t=1$ then by changing variable $t=1-s$, these boundary conditions can be transformed to the boundary conditions of the form \eqref{bc1}. One set of boundary conditions among the above mentioned conditions is $u'(0)=u(1)=u'(1)=0$ considered in \cite{Feng09}. Therefore, from now on we shall consider only the boundary conditions of the form \eqref{bc1}.\par
To investigate the BVPs \eqref{eqfull}, \eqref{bc1} we use a novel approach based on the reduction of them to operator equations for the nonlinear terms but not for the functions to be sought. This approach was used by ourselves to some boundary value problems for fourth order nonlinear equations in very recent works \cite{A-Quy,A-L-Q}. Here,
by this approach we have established the existence, uniqueness, positivity and monotony of solutions and the convergence of the iterative method for approximating the solutions of the problems \eqref{eqfull}-\eqref{bc1}  under some easily verified conditions in bounded domains. These conditions are much simpler and weaker than those of other authors for studying solvability of particular cases of the problems before by using different methods. Many examples illustrate the obtained theoretical results.

\section{Existence results }
For convenience we rewrite the problem \eqref{eqfull}-\eqref{bc1} in the form 
\begin{equation}\label{eq1}
\begin{split}
u'''(t)&=f(t,u(t),u'(t),u''(t)), \quad 0<t<1\\
B_1[u]&=B_2[u]=B_3[u]=0,
\end{split}
\end{equation}
where $B_1[u], B_2[u], B_3[u]$ are defined by \eqref{bc1}. We shall associate this problem with an operator equation as follows.\\
For functions $\varphi (x) \in C[0, 1]$ consider the nonlinear operator $A$ defined by
\begin{equation}\label{defA}
(A\varphi )(t)=f(t,u(t), u'(t), u''(t)),
\end{equation}
where $u(t)$ is a solution of the problem
\begin{equation}\label{eq2}
\begin{split}
u'''(t)&=\varphi (t), \quad 0<t<1\\
B_1[u]&=B_2[u]=B_3[u]=0.
\end{split}
\end{equation}
It is easy to verify the following
\begin{prop}\label{prop1}   If the function $\varphi (x)$ is a fixed point of the operator $A$, i.e., $\varphi (t)$ is a solution of the operator equation
\begin{equation}\label{opereq}
\varphi = A\varphi , 
\end{equation}
then the function $u(t)$ determined from the boundary value problem \eqref{eq2} solves the problem \eqref{eq1}. Conversely, if $u(t)$ is a solution of the boundary value problem \eqref{eq1} then the function
\begin{equation*}
\varphi(t)=f(t,u(t), u'(t), u''(t))
\end{equation*}
is a fixed point of the operator $A$ defined above by \eqref{defA}, \eqref{eq2}.
\end{prop}
Thus, the solution of the original problem \eqref{eq1} is reduced to the solution of the operator equation \eqref{opereq}.\par
Now consider the problem \eqref{eq2}. Suppose that the Green function of it exists and is denoted by $G(t,s)$. Then the unique solution of the problem is represented in the form
\begin{equation}\label{equ}
u(t)=\int _0 ^1 G(t,s)\varphi (s) ds.
\end{equation}
By differentiation of both sides of the above formula we obtain
\begin{equation}\label{equ'u''}
u'(t)= \int _0 ^1 G_1(t,s)\varphi (s) ds,\quad  u''(t)= \int _0 ^1 G_2(t,s)\varphi (s) ds,
\end{equation}
where $G_1(t,s)=G'_{t}(t,s)$ is a function continuous in the square $Q= [0,1]^2$ and $G_2(t,s)=G''_{tt}(t,s)$ is continuous in the square $Q$ except for the line $t=s$.\\
Further, let
\begin{equation}\label{valGreen}
\begin{split}
\max _{0\le t\le 1}\int _0 ^1 |G(t,s)| ds &=M_0 \\
\max _{0\le t\le 1}\int _0 ^1 |G_1(t,s)| ds &=M_1, \;
\max _{0\le t\le 1}\int _0 ^1 |G_2(t,s)| ds =M_2.
\end{split}
\end{equation}
Next, for each fixed real number $M>0$ introduce the domain
\begin{equation}\label{defDM}
\mathcal{D}_M=\{ (t,x,y,z)| \ 0\leq t\leq 1, \,\,|x| \leq M_0M, \,\, |y| \leq M_1M, \,\, 
 |z| \leq M_2M \},
\end{equation}
and as usual, by $B[O,M]$ we denote the closed ball of radius $M$ centered at $0$ in the space of continuous in $[0, 1]$ functions, namely,
\begin{equation*}
B[O,M]=\{ \varphi \in C[0,1]| \ \| \varphi \| \leq M \},
\end{equation*}
where
\[\| \varphi \|=  \max_{0 \leq t \leq 1} |\varphi (t)|. \]
\begin{thm}\label{theorem1}
Suppose that there exists a number $M>0$ such that the function $f(t,x,y,z)$ is continuous and bounded by $M$ in the domain $\mathcal{D}_M$, i.e., 
\begin{equation}\label{cond1}
|f(t,x,y,z)| \leq M,
\end{equation}
for any $(t,x,y,z) \in \mathcal{D}_M .$\par
Then, the problem \eqref{eq1} has a solution $u(t)$ satisfying 
\begin{equation}\label{estimate1}
|u(t)| \leq M_0M, \; |u'(t)| \leq M_1M,\; |u''(t)| \leq M_2M \text{ for any } 0 \le t \le 1.
\end{equation}
\end{thm}
\noindent {\bf Proof.}
Having in mind Proposition \ref{prop1} the theorem will be proved if we show that the operator $A$ associated with the problem \eqref{eq1} has a fixed point. For this purpose, it is not difficult to show that the operator $A$ 
maps the closed ball $B[0, M]$ into itself. Next, from the compactness of integral operators \eqref{equ}, \eqref{equ'u''}, which put each $\varphi \in C[0, 1]$ in correspondence to the functions $u, u', u''$, respectively
 \cite[Sec. 31]{Kolm} (see APPENDIX) and the continuity of the function $f(t,x,y,z)$ it follows that $A$ is a compact operator in the Banach space $C[0,1]$. By the Schauder Fixed Point Theorem \cite{Zeidler} the operator $A$ has a fixed point in $B[0, M]$. The estimates \eqref{estimate1} hold due to the equalities \eqref{equ}, \eqref{equ'u''} and \eqref{valGreen}. 
The theorem is proved. \quad $\square$ \par
\medskip
Now suppose that the Green function $G(x,t)$ and its first derivative $G_1(x,t)$ are of constant signs in the square $Q= [0, 1]^2$. Let's adopt the following convention:\\
For a function $H(x,t)$ defined  and having a constant sign in the square $Q$ we define
$$
\sigma (H)=sign (H(t,s)) = \left\{
\begin{array}{ll}
1, \quad \text{ if } H(t,s) \ge 0 ,\\
-1, \quad \text{ if } H(t,s) \le 0.
\end{array}
\right.
$$

In order to investigate the existence of positive solutions of the problem \eqref{eqfull},\eqref{bc1} we introduce the notations
\begin{equation}\label{eqD+}
\begin{split}
\mathcal{D}_M^+=\{ (t,x,y,z)| \ &0\leq t\leq 1, \,\,0 \leq x \leq M_0M, \,\,\\
 &0 \leq \sigma (G)\sigma (G_1)y\le M_1M, \,\,  |z| \leq M_2M \}
\end{split}
\end{equation}
and 
\begin{equation}\label{eqS}
S_M=\{ \varphi \in C[0,1]| \ 0\le  \sigma (G)\varphi  \leq M \}.
\end{equation}
\begin{thm} \label{theorem2}(Existence of constant sign solution)
Suppose that there exists a number $M>0$ such that the function $f(t,x,y,z)$ is continuous and 
\begin{equation}\label{cond1a}
0\le \sigma (G)f(t,x,y,z) \leq M
\end{equation}
for any $(t,x,y,z) \in \mathcal{D}_M^+ $.
Then, the problem \eqref{eqfull},\eqref{bc1} has a monotone nonnegative solution $u(t)$ satisfying
\begin{equation}\label{estimate2}
\begin{split}
0\le  u(t) \leq M_0M, \,\, 0\le \sigma (G)\sigma (G_1)u'(t) \leq M_1M, \,\, 
 |u''(t)| \leq M_2M.
\end{split}
\end{equation}

In addition, if $\sigma (G) \sigma (G_1) =1 $ then the problem has a nonnegative, increasing solution, and if $\sigma (G) \sigma (G_1) =-1 $ then the problem has a nonnegative, decreasing solution.
\end{thm} 
\noindent {\bf Proof.} The proof of the existence of monotone nonnegative solution of the problem is similar to that of solution in Theorem \ref{theorem1} with the replacements of $\mathcal{D_M}$ by $\mathcal{D}_M^+$, $B[0, M]$ by $S_M$ and the condition \eqref{cond1} by the condition \eqref{cond1a}. 
From the estimates \eqref{estimate2} it is obvious that if $\sigma (G) \sigma (G_1) =1 $ then $u'(t) \ge 0$, consequently, the solution is increasing function, otherwise, if $\sigma (G) \sigma (G_1) =-1 $ then $u'(t) \le 0$, therefore, the solution is decreasing function. The theorem is proved.

\begin{thm}\label{theorem3}(Existence and uniqueness of solution)
Assume that there exist numbers
  $M,L_0, L_1, L_2 \geq 0$ such that
\begin{equation}\label{cond1'}
|f(t,x,y,z)| \leq M,
\end{equation}
\begin{multline}\label{Lips}
|f(t,x_2,y_2,z_2)-f(t,x_1,y_1,z_1)| \leq \\ L_0|x_2-x_1|+ L_1|y_2-y_1|+L_2|z_2-z_1|       
\end{multline}
for any $(t,x,y,z), (t,x_i,y_i,z_i) \in \mathcal{D}_M \ (i=1,2)$ and
\begin{equation}\label{eqq}
q:=L_0M_0+ L_1M_1+L_2M_2<1.
\end{equation}
Then, the problem \eqref{eqfull},\eqref{bc1} has a unique solution $u(t)$ such that $|u(t)| \leq M_0M,$ $|u'(t)| \leq M_1M, \,\, |u''(t)| \leq M_2M $ for any $0 \le t \le 1$.
\end{thm}
\noindent {\bf Proof.} It is easy to show that under the conditions of the theorem, the operator $A$ associated with the problem \eqref{eqfull},\eqref{bc1} is a contraction mapping from the closed ball $B[0, M]$ into itself. 
By the contraction principle the operator $A$ has a unique fixed point in $B[O,M]$, which corresponds to a unique solution $u(t)$ of the problem \eqref{eqfull},\eqref{bc1}.\\
The estimates for $u(t)$ and its derivatives are obtained as in Theorem \ref{theorem1}. Thus, the theorem is proved. \quad $\square$\par
\medskip
Analogously, we have the following theorem for the existence and uniqueness of constant sign solution of the problem \eqref{eqfull},\eqref{bc1}.
\begin{thm}\label{theorem4}(Existence and uniqueness of constant sign solution)
Assume that all the conditions of Theorem \ref{theorem2} are satisfied in the domain
$\mathcal{D}_M^+ $. Moreover, assume that there exist numbers
  $L_0, L_1, L_2 \geq 0$ such that the function $f(t,x,y,z)$ satisfies the Lipschitz conditions \eqref{Lips}, \eqref{eqq}.
Then, the problem \eqref{eqfull},\eqref{bc1} has a unique monotone nonnegative  solution $u(t)$ satisfying \eqref{estimate2}.
\end{thm}
\section{Iterative method }
Consider the following iterative method for solving the problem \eqref{eqfull}, \eqref{bc1}:
\begin{enumerate}
\item Given a starting approximation $\varphi _0 \in B[0, M]$, say 
\begin{equation}\label{iter1}
\varphi_0(t)=f(t,0,0,0).
\end{equation}
\item Knowing $\varphi_k \; (k=0,1,...)$ compute
\begin{equation}\label{iter2}
\begin{split}
u_k(t)&= \int _0 ^1 G(t,s)\varphi_k(s) \ ds,\\
y_k(t)&= \int _0 ^1 G_1(t,s)\varphi_k(s) \ ds,\\
z_k(t)&= \int _0 ^1 G_2(t,s)\varphi_k(s) \ ds.
\end{split}
\end{equation}
\item Update the new approximation
\begin{equation}\label{iter3}
\varphi_{k+1}(t)=f(t,u_k(t),y_k(t),z_k(t)).
\end{equation}
\end{enumerate}
The above iterative process indeed is the successive approximation of the fixed point of the operator $A$ associated with the problem  \eqref{eqfull},\eqref{bc1}. Therefore, it converges with the rate of geometric progression and there is the estimate
\begin{equation} \label{iter4}
\|\varphi _k - \varphi \| \leq p_k,
\end{equation}
where $\varphi $ is the fixed point of $A$ and
\begin{equation*}
p_k=\dfrac{q^k}{1-q}\| \varphi _1 -\varphi _0\|.
\end{equation*}
Taking into account the representations \eqref{equ}, \eqref{equ'u''} and \eqref{iter2}, from the above estimate we obtain the following error estimates for the approximate solution $u_k$ and its derivatives
\begin{equation}\label{iter5}
\|u_k-u\| \leq M_0p_k, \quad  \|u'_k-u'\| \leq M_1p_k,\quad
\|u''_k-u''\| \leq M_2p_k,
\end{equation}
where $u$ is the exact solution of the problem \eqref{eqfull}, \eqref{bc1}.

\section{Some particular cases and examples}
Consider some particular cases of the general boundary value problem \eqref{eqfull},\eqref{bc1}, which cover the problems studied by other authors using different methods. For each case, the theoretical results obtained in the previous section will be illustrated on examples, including some examples considered before by other authors. In numerical realization of the proposed iterative method,  for computing definite integrals the trapezium formula with second order accuracy is used.
In all examples, numerical computations are performed on the uniform grid on the interval $[0, 1]$ with the gridsize $h=0.01$ until $\| \varphi_ k - \varphi _{k-1}\| \le 10^{-6}$. The number of iterations for reaching the above accuracy will be indicated.

 \par
\subsection{Case 1.} Consider the problem
\begin{equation}\label{case1}
\begin{split}
u^{(3)}(t) &=f(t, u(t), u'(t), u''(t)), \quad 0 < t < 1, \\
u(0)&=u'(0)=u'(1)=0.
\end{split}
\end{equation}
The Green function associated with the above problem has the form
\begin{equation}\label{green1}
\begin{aligned}
G(t,s)=\left\{\begin{array}{ll}
\dfrac{s}{2}(t^2-2t+s), \quad 0\le s \le t \le 1,\\
\, \, \dfrac{t^2}{2}(s-1), \quad 0\le t \le s \le 1.\\
\end{array}\right.
\end{aligned}
\end{equation}
After differentiation of $G(t,s)$ we obtain
\begin{equation*}
G_1(t,s)=\left\{\begin{array}{ll}
s(t-1), \quad 0\le s \le t \le 1,\\
t(s-1), \quad 0\le t \le s \le 1,\\
\end{array}\right.
\end{equation*}
\begin{equation*}
G_2(t,s)=\left\{\begin{array}{ll}
s,& \quad 0\le s \le t \le 1,\\
s-1,& \quad 0\le t \le s \le 1.\\
\end{array}\right.
\end{equation*}
It is obvious that
\begin{equation*}
G(t,s) \le 0, \; G_1(t,s) \le 0, \; G_2(t,s) \ge 0, \; 0 \le t, s \le 1
\end{equation*}
and we have
\begin{equation*}
\begin{aligned}
&M_0= \max _{0\le t\le 1} \int _0 ^1 |G(t,s)| \ ds =\dfrac{1}{12},\quad M_1= \max _{0\le t\le 1} \int _0 ^1 |G_1(t,s)| \ ds =\dfrac{1}{8},\\
&M_2= \max _{0\le t\le 1} \int _0 ^1  G_2(t,s)  \ ds =\dfrac{1}{2}.\\
\end{aligned}
\end{equation*}
\noindent {\bf Example 4.1.1 } (Example 7 in \cite{Yao-Feng})
Consider the problem
\begin{equation}\label{Yao-Feng7}
\begin{split}
u^{(3)}(t) &=-e^{u(t)}, \quad 0 < t < 1, \\
u(0)&=u'(0)=u'(1)=0.
\end{split}
\end{equation}
Yao and Feng \cite{Yao-Feng} using the lower and upper solutions method and the fixed point theorem on cones  proved that the above problem has a solution $u(t)$ such that $\| u\| \le 1, \, u(t) >0$ for $t\in (0, 1)$ and $u(t)$ is an increasing function. Here, using the theoretical results obtained in the previous section we establish the results which are more strong than the above results.\par
Indeed, for the problem \eqref{Yao-Feng7} $f=f(t,x)=-e^x$. In the domain 
\begin{equation*}
\begin{split}
\mathcal{D}_M^+=\Big \{ (t,x)| \ 0\leq t\leq 1, \,\, 
0\le  x \leq \frac{M}{12}
  \Big \}
\end{split}
\end{equation*}
there hold $ -e^{M/12} \le f(t,x) \le 0$. So, with the choice $M=1.1$ we have $-M \le f(t,x) \le 0$. Further, in $\mathcal{D}_M^+$ the function $f(t,x)$ satisfies the Lipschitz condition with $L_0=e^{M/12}=1.096$. Therefore, $q=L_0/12=0.0913$. By Theorem \ref{theorem4} the problem has a \emph{unique} monotone positive solution $u(t)$ satisfying the estimates 
\begin{equation*}
\begin{split}
0&\le u(t) \le \dfrac{M}{12}=\dfrac{1.1}{12}=0.0917, \, 0\le u'(t) \le \dfrac{M}{8}=\dfrac{1.1}{8}=0.1357, \,\\
&|u''(t)| \le \dfrac{M}{2}=\dfrac{1.1}{2}=0.55.
\end{split}
\end{equation*}
Clearly, these results are better than those in \cite{Yao-Feng}.\par
\begin{figure}
\centering
\includegraphics[width=0.8\textwidth]{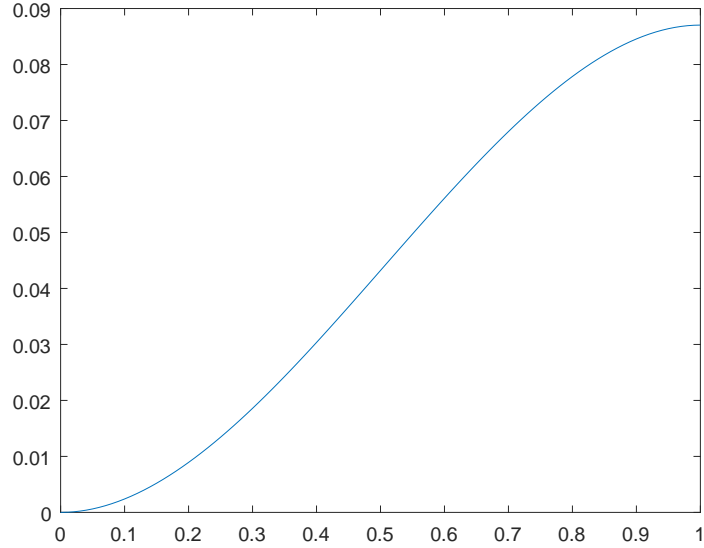}
\caption{The graph of the approximate solution in Example 4.1.1.} \label {Fig1}
\end{figure}

The numerical solution of the problem obtained by the iterative method \eqref{iter1}-\eqref{iter3} after 5 iterations is depicted in Figure \ref{Fig1}.
From this figure \ref{Fig1} it is clear that the solution is monotone, positive and is bounded  by $0.0917$ as shown above by the theory.\\

\begin{figure}
\centering
\includegraphics[width=0.8\textwidth]{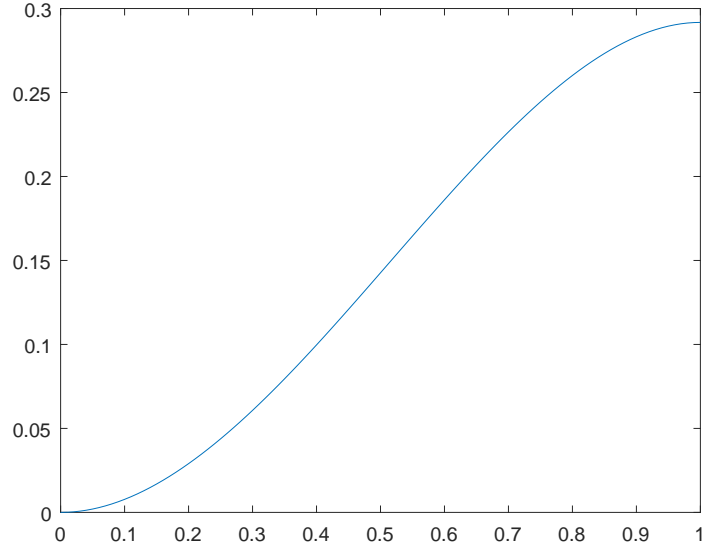}
\caption{The graph of the approximate solution in Example 4.1.2.} \label {Fig2}
\end{figure}

\noindent {\bf Example 4.1.2 } (Example 8 in \cite{Yao-Feng})
Consider the problem
\begin{equation}\label{Yao-Feng8}
\begin{split}
u^{(3)}(t) &=-\frac{5u^3(t)+4u(t)+3}{u^2(t)+1}, \quad 0 < t < 1, \\
u(0)&=u'(0)=u'(1)=0.
\end{split}
\end{equation}
Yao and Feng in \cite{Yao-Feng} showed that the above problem has a solution $u(t)$ such that $ u(t) >0$ for $t\in (0, 1)$ and $u(t)$ is an increasing function. Similarly as in Example 4.1.1 we established that the problem \eqref{Yao-Feng8} has a \emph{unique} monotone positive solution $u(t)$ satisfying
$$ 0\le u(t) \le 0.3417, \, 0\le u'(t) \le 0.5125, \, |u''(t)| \le 2.05.$$
The numerical solution of the problem obtained by the iterative method \eqref{iter1}-\eqref{iter3} after 8 iterations is depicted in Figure \ref{Fig2}.
From this figure \ref{Fig2} it is clear that the solution is monotone, positive and is bounded  by $0.3417$ as shown above by the theory.\\


\noindent {\bf Example 4.1.3 } (Example 4.2 in \cite{Feng-Liu})
Consider the problem
\begin{equation}\label{Feng-Liu4.2}
\begin{split}
u^{(3)}(t) &=-e^{u(t)} -e^{u'(t)}, \quad 0 < t < 1, \\
u(0)&=u'(0)=u'(1)=0.
\end{split}
\end{equation}
Using the lower and upper solutions method and a new maximum principle, Feng and Liu in \cite{Feng-Liu}  established that the above problem has a solution $u(t)$ such that $\| u\| \le 1, \, u(t) >0$ for $t\in (0, 1)$ and $u(t)$ is an increasing function. Here, using Theorem \ref{theorem4} with the choice $M=2.7$ we conclude that the problem has a \emph{unique} monotone positive solution $u(t)$ satisfying the estimates 
$$  0\le u(t) \le 0.2250, \, 0\le u'(t) \le 0.3375, \, |u''(t)| \le 1.350$$

\begin{figure}
\centering
\includegraphics[width=0.8\textwidth]{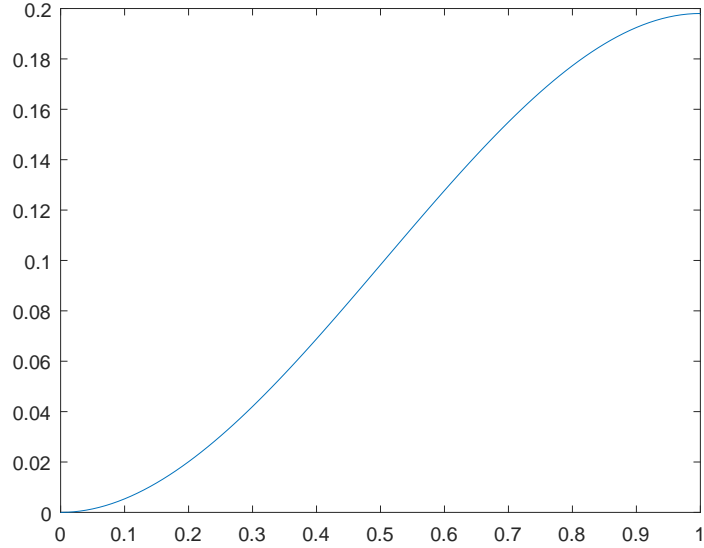}
\caption{The graph of the approximate solution in Example 4.1.3.} \label {Fig3}
\end{figure}

The numerical solution of the problem obtained by the iterative method \eqref{iter1}-\eqref{iter3} after 9 iterations is depicted in Figure \ref{Fig3}.
From this figure \ref{Fig3} it is clear that the solution is monotone, positive and is bounded  by $0.2250$ as shown above by the theory.\\

\begin{rem}
 It should be emphasized that in \cite{Yao-Feng} and \cite{Feng-Liu} the authors used one very important assumption, which means that the nonlinear functions $f(t,x)$ or $f(t,x,y)$  satisfy one-side Lipschitz condition in $x$ or $x, y$ in the whole space $\mathbb{R}$ or $\mathbb{R}^2$, respectively. If now change the sign of the right-hand sides then this condition is not satisfied. Therefore, it is impossible to say anything about the solution of the problem.
 But Theorem \ref{theorem3} ensures the existence and uniqueness of a solution. Moreover, in a similar way as in Theorem \ref{theorem3} it is possible conclude that this solution is nonpositive.
\end{rem}

\subsection{Case 2.} Consider the problem
\begin{equation}\label{case2}
\begin{split}
u^{(3)}(t) &=f(t, u(t), u'(t), u''(t)), \quad 0 < t < 1, \\
u(0)&=u'(0)=u''(1)=0.
\end{split}
\end{equation}
In \cite{Hopkins} under the assumptions that the function $f(t,x,y,z)$ defined on $[0, 1] \times \mathbb{R}^3 \rightarrow \mathbb{R}$ is $L_p$-Caratheodory, and there exist functions $\alpha , \beta , \gamma , \delta \in L_p[0, 1], \, p \ge 1$, such that
\begin{equation*}
|f(t,x,y,z)\le \alpha (t) x +\beta (t) y + \gamma (t) z + \delta (t)|, \quad t\in (0, 1)
\end{equation*}
and 
\begin{equation*}
A_0 \| \alpha\|_p + A_1 \| \beta\|_p + \| \gamma\|_p <1,
\end{equation*}
where $A_0, A_1$ are some constants depending on $p$, the problem has at least one solution. The tool used is the Leray-Schauder continuation principle. No examples are given for illustrating the theoretical results.\par
Here, assuming that the function $f(t,x,y,z)$ is continuous, we establish the existence of unique solution by Theorem \ref{theorem4}.
For the problem \eqref{case2} the Green function is
\begin{equation}\label{green2}
\begin{aligned}
G(t,s)=\left\{\begin{array}{ll}
-st +\dfrac{s^2}{2}, \quad 0\le s \le t \le 1,\\
\, \, -\dfrac{t^2}{2}, \quad 0\le t \le s \le 1.\\
\end{array}\right.
\end{aligned}
\end{equation}
The first and the second derivatives of this function are
\begin{equation*}
G_1(t,s)=\left\{\begin{array}{ll}
-s, \quad 0\le s \le t \le 1,\\
-t, \quad 0\le t \le s \le 1,\\
\end{array}\right.
\end{equation*}
\begin{equation*}
G_2(t,s)=\left\{\begin{array}{ll}
0,& \quad 0\le s \le t \le 1,\\
-1,& \quad 0\le t \le s \le 1.\\
\end{array}\right.
\end{equation*}
It is easy to see that
\begin{equation*}
G(t,s) \le 0, \; G_1(t,s) \le 0, \; 0 \le t, s \le 1
\end{equation*}
and
\begin{equation*}
\begin{aligned}
&M_0= \max _{0\le t\le 1} \int _0 ^1 |G(t,s)| \ ds =\dfrac{1}{3},\quad M_1= \max _{0\le t\le 1} \int _0 ^1 |G_1(t,s)| \ ds =\dfrac{1}{2},\\
&M_2= \max _{0\le t\le 1} \int _0 ^1  |G_2(t,s)|  \ ds =1.
\end{aligned}
\end{equation*}

\noindent {\bf Example 4.2.1 }
Consider the following problem
\begin{equation}\label{dqa1}
\begin{aligned}
u'''(t)&=-\dfrac{1}{36} \big ( u'(t)  \big )^2+\dfrac{1}{24}u(t) u''(t) +\dfrac{1}{4}t^2-6, \quad 0\le t \le 1,\\
u(0)&=u'(0)=u''(1)=0.
\end{aligned}
\end{equation}
In this example 
\begin{equation*}
\begin{aligned}
f(t,x,y,z)&=-\dfrac{1}{36}y^2 +\dfrac{1}{24}xz +\dfrac{1}{4}t^2-6.\\
\end{aligned}
\end{equation*}
It is possible to verify that with $M=7.5,\; L_1=0.3125, \; L_2=0.2083,\; L_3=0.1042$, the conditions of Theorem \ref{theorem4} are met, therefore, the problem \eqref{dqa1} has a unique positive solution satisfying the estimates $0 \le u(t) \le 2.5, \; 0 \le u'(t) \le 3.75 , \; |u''(t)| \le 7.5$. \par

\begin{figure}
\centering
\includegraphics[width=0.8\textwidth]{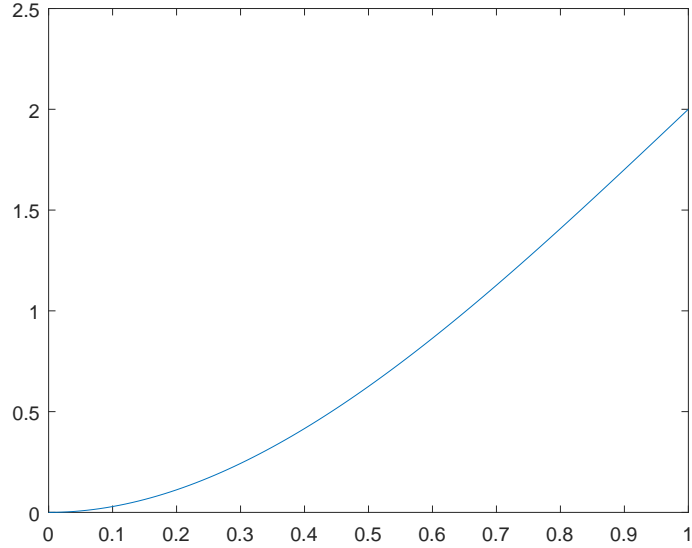}
\caption{The graph of the approximate solution in Example 4.2.1.} \label {Fig4.2.1}
\end{figure}

The numerical solution of the problem obtained by the iterative method \eqref{iter1}-\eqref{iter3} after $5$ iterations is depicted in Figure \ref{Fig4.2.1}.
From this figure it is clear that the solution is bounded  by $2.5$ as shown above by the theory.\par
It is interesting that the problem \eqref{dqa1} has the  exact solution $u(t)=-t^3+3t^2$.
This solution satisfies the exact estimates
$0 \le u(t) \le 2, \; 0 \le u'(t) \le 3 , \; 0 \le u''(t) \le 6$ for $0 \le t \le 1$, which are better than the theoretical estimates above.
On the grid with the gridsize $h=0.01$ the maximal deviation of the obtained approximate solution and the exact solution is $3.7665e-04$. 
\subsection{Case 3.} Consider the problem
\begin{equation}\label{case3}
\begin{split}
u^{(3)}(t) &=f(t, u(t), u'(t), u''(t)), \quad 0 < t < 1, \\
u(0)&=u'(1)=u''(1)=0.
\end{split}
\end{equation}
Under the conditions similar to those in the previous case, Hopkins and Kosmatove  in \cite{Hopkins} established the existence of a solution of the problem without illustrative examples. Very recently, in \cite{Li-Li} Li Yongxiang and  Li Yanhong studied the existence of positive solutions of the problem \eqref{case3} under conditions on the growth of the function $f(t,x,y,z)$ as $|x|+|y|+|z|$ tends to zero and infinity, including a Nagumo-type condition on $y$ and $z$. The tool used is the fixed point index theory on cones.
\par
Here, assuming that the function $f(t,x,y,z)$ is continuous, we can establish the existence results by the theorems in the Section 2.
For the problem \eqref{case3} the Green function is
\begin{equation}\label{green3}
\begin{aligned}
G(t,s)=\left\{\begin{array}{ll}
\, \, \dfrac{s^2}{2}, \quad 0\le s \le t \le 1,\\
st -\dfrac{t^2}{2}, \quad 0\le t \le s \le 1.\\
\end{array}\right.
\end{aligned}
\end{equation}
The first and the second derivatives of this function are
\begin{equation*}
G_1(t,s)=\left\{\begin{array}{ll}
0, \quad 0\le s \le t \le 1,\\
s-t, \quad 0\le t \le s \le 1,\\
\end{array}\right.
\end{equation*}
\begin{equation*}
G_2(t,s)=\left\{\begin{array}{ll}
0,& \quad 0\le s \le t \le 1,\\
-1,& \quad 0\le t \le s \le 1.\\
\end{array}\right.
\end{equation*}
It is easy to see that
\begin{equation*}
G(t,s) \ge 0, \; G_1(t,s) \ge 0, \; 0 \le t, s \le 1
\end{equation*}
and to obtain
\begin{equation*}
M_0=\dfrac{1}{6}, \, M_1= \dfrac{1}{2}, \, M_2 =1.
\end{equation*}

\noindent {\bf Example 4.3.1 }
Consider the following problem
\begin{equation}\label{dqa}
\begin{aligned}
u'''(t)&=\dfrac{1}{18} \big ( u'(t)  \big )^2-\dfrac{1}{12}u(t) u''(t) +\dfrac{1}{2}t+\dfrac{11}{2}, \quad 0\le t \le 1,\\
u(0)&=u'(1)=u''(1)=0.
\end{aligned}
\end{equation}
In this example 
\begin{equation*}
\begin{aligned}
f(t,x,y,z)&=\dfrac{1}{18}y^2 -\dfrac{1}{12}xz +\dfrac{1}{2}t+\dfrac{11}{2}.\\
\end{aligned}
\end{equation*}
It is possible to verify that with $M=8,\; L_1=\dfrac{2}{3}, \; L_2=\dfrac{4}{9},\; L_3=\dfrac{1}{9}$ all the conditions of Theorem \ref{theorem4} are met, therefore, the problem \eqref{dqa} has a unique positive solution, which is increasing and satisfies the estimates $0 \le u(t) \le \dfrac{4}{3}, \; 0 \le u'(t) \le 4 , \; -8 \le u''(t) \le 0$. \par
\begin{figure}
\centering
\includegraphics[width=0.8\textwidth]{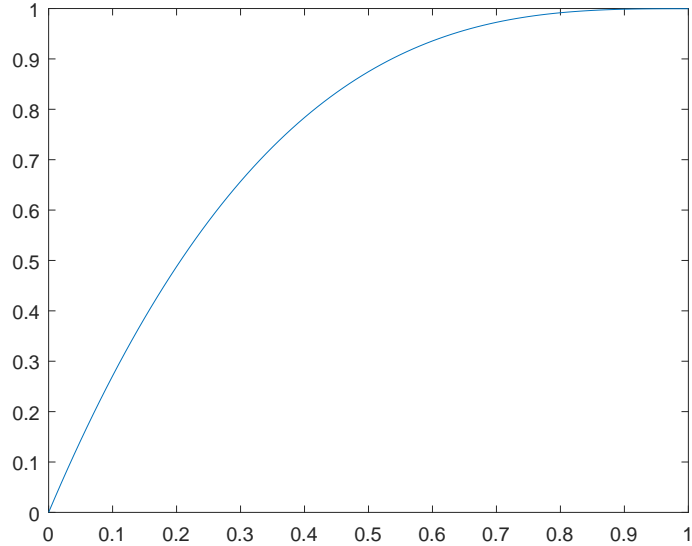}
\caption{The graph of the approximate solution in Example 4.3.1.} \label {Fig4}
\end{figure}

The numerical solution of the problem obtained by the iterative method \eqref{iter1}-\eqref{iter3} after $6$ iterations is depicted in Figure \ref{Fig4}.
From this figure \ref{Fig4} it is clear that the solution is monotone, positive and is bounded  by $4/3$ as shown above by the theory.\par

It is possible to verify that the function $u(t)=t^3-3t^2+3t$ is the  exact solution of the problem \eqref{dqa}. This solution is positive, increasing and  satisfies the exact estimates
$0 \le u(t) \le 1, \; 0 \le u'(t) \le 3 , \; -6 \le u''(t) \le 0$ for $0 \le t \le 1$, which are better than the theoretical estimates above.
On the grid with the gridsize $h=0.01$ the maximal deviation of the obtained approximate solution and the exact solution is $3.6256e-04$.

%

\subsection{Case 4.} Consider the problem
\begin{equation}\label{case4}
\begin{split}
u^{(3)}(t) &=f(t, u(t), u'(t), u''(t)), \quad 0 < t < 1, \\
u(0)&=u''(0)=u'(1)=0.
\end{split}
\end{equation}
Using the lower and upper solutions method and Schauder fixed theorem on cones, Bai \cite{Bai} established the existence of a solution under complicated conditions on the right-hand side function.\par

For the problem \eqref{case4} the Green function is
\begin{equation}\label{green4}
\begin{aligned}
G(t,s)=\left\{\begin{array}{ll}
\dfrac{t^2}{2}-t+ \dfrac{s^2}{2}, \quad 0\le s \le t \le 1,\\
t(s-1), \quad 0\le t \le s \le 1.\\
\end{array}\right.
\end{aligned}
\end{equation}
The first and the second derivatives of this function are
\begin{equation*}
G_1(t,s)=\left\{\begin{array}{ll}
t-1, \quad 0\le s \le t \le 1,\\
s-1, \quad 0\le t \le s \le 1,\\
\end{array}\right.
\end{equation*}
\begin{equation*}
G_2(t,s)=\left\{\begin{array}{ll}
1,& \quad 0\le s \le t \le 1,\\
0,& \quad 0\le t \le s \le 1.\\
\end{array}\right.
\end{equation*}
Obviously,
\begin{equation*}
G(t,s) \le 0, \; G_1(t,s) \le 0, \; 0 \le t, s \le 1
\end{equation*}
and it is easy to obtain
\begin{equation*}
M_0=\dfrac{1}{3}, \, M_1= \dfrac{1}{2}, \, M_2 =1.
\end{equation*}
In view of the above facts concerning the Green function, using theorems in the previous section we can establish the results on the existence of solution of the problem \eqref{case4}. \\

\begin{figure}
\centering
\includegraphics[width=0.8\textwidth]{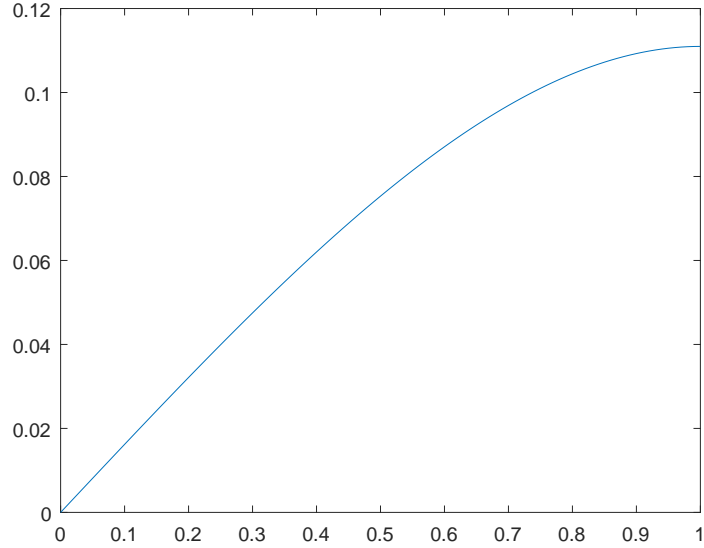}
\caption{The graph of the approximate solution in Example 4.4.1.} \label {Fig5}
\end{figure}

\noindent {\bf Example 4.4.1 } (Example 3.5 in \cite{Bai})
\begin{equation}\label{Bai3.5}
\begin{split}
u^{(3)}(t) &=-\dfrac{1}{4} \big [ t+e^{u(t)}+(u'(t))^2 +u''(t)   \big ]  , \quad 0 < t < 1, \\
u(0)&=u''(0)=u'(1)=0.
\end{split}
\end{equation}
Defining
\begin{equation*}
\begin{split}
\mathcal{D}_M^+=\{ (t,x,y,z)| \ 0\leq t\leq 1, \,\,0 \leq x \leq \frac{M}{3}, \,
 0 \leq y\le\dfrac{M}{2}, \,\,  |z| \leq M \},
\end{split}
\end{equation*}
for $M=0.835$ we have
\begin{equation*}
-M \le f(t,x,y,z)= -\dfrac{1}{4} \big [ t+e^x+y^2 +z   \big ]\le 0
\end{equation*}
Further, it is easy to calculate the Lipschitz coefficients of $f(t,x,y,z)$: 
$$ L_0=\dfrac{1}{4}e^{M/3}=0.3302,\,   L_1= \dfrac{M}{4}=0.2087,\,  L_2=1.$$
Therefore, $q=L_0/3 +L_1/2+L_2=0.4851<1.$ By Theorem \ref{theorem4} the problem has a unique monotone positive solution $u(t)$ such that
\begin{equation*}
\begin{split}
0\le u(t) \le M/3=0.2783, \, 0\le u'(t) \le M/2=0.5,\, |u''(t)|\le 1. 
\end{split}
\end{equation*}
Notice that in \cite{Bai} Bai could only conclude that the problem \eqref{Bai3.5} has a positive solution.\\
The numerical solution of the problem obtained by the iterative method \eqref{iter1}-\eqref{iter3} after $5$ iterations is depicted in Figure \ref{Fig5}.
From this figure \ref{Fig5} it is clear that the solution is monotone, positive and is bounded  by $0.2783$ as shown above by the theory.\\

\section{Conclusion}
In this paper we have proposed a novel approach to study fully third order differential equation with general two-point linear boundary conditions. The approach is based on the reduction of boundary value problems to fixed point problems for nonlinear operators for the right-hand sides of the equation but not for the function to be sought. The results are that we have established the existence, uniqueness, positivity and monotony of solution under the conditions, which are simpler and easier to verify than those of other authors. The applicability and advantages of the proposed approach are illustrated on some examples taken from the papers of other authors, where our approach gives better results. 
\section*{Acknowledgments}
This work is supported by Vietnam National Foundation for Science and Technology
Development (NAFOSTED) under the grant  number 102.01-2017.306. \\

\section*{APPENDIX}
 In the space $C[a, b]$ consider the operator $y=Ax$ defined by the formula
\begin{equation*}
y(t)=\int _a^b K(t,s) x(s) ds.
\end{equation*}
\begin{thm}(see \cite[Sec. 31]{Kolm})
The above formula  defines a compact operator in the
space $C[a, b]$ if the function $K(t, s)$ is bounded on the square $a \le t \le b, 
a \le s \le b$ and all points of discontinuity of the function $K(s, t)$ lie on a finite
number of curves
$$s = \varphi_k(t),\quad k = 1, 2. ... , n,$$
where $\varphi_k(t)$ the are continuous functions.
\end{thm}

\end{document}